\documentclass[12pt]{article}
\usepackage[psamsfonts]{amssymb}

\makeatletter
\def\@cite#1#2{$^{\mbox{\scriptsize #1\if@tempswa , #2\fi}}$}
\makeatother

\def\thebibliography#1{\section*{References\markboth
  {REFERENCES}{REFERENCES}}\list
  {$^{\mbox{\scriptsize\arabic{enumi}}}$\hspace{-2mm}}{\settowidth\labelwidth{[#1]}\leftmargin\labelwidth    
  \advance\leftmargin\labelsep 
  \usecounter{enumi}}
  \def\newblock{\hskip .11em plus .33em minus -.07em}
  \sloppy
  \sfcode`\.=1000\relax}
\def\t#1{\tilde #1}
\def\p{\partial}
\def\h#1{\hat #1}

\def\kbar{{\mathchar'26\mkern-9muk}}  
\def\bra#1{\langle #1 \vert}
\def\ket#1{\vert #1 \rangle}

\def\b#1{{\mathbb #1}}
\def\c#1{{\cal #1}}

\def\Dirac{{\raise0.09em\hbox{/}}\kern-0.69em D}
\def\sgn{\mbox{sgn}}
\def\sq{\hbox{\rlap{$\sqcap$}$\sqcup$}}
\def\etal{{\it et al.}\ }
\newcounter{eg}
\newcounter{chapter}
\newtheorem{eg}{Example}[chapter]    
\def\beg{\begin{eg}\rm}
\def\eeg{\hfill\sq\end{eg}}
\def\exterior{{{\raise0.2em\hbox{$\scriptstyle\bigwedge$}}{}}}
\def\k{\kern-.1em\mathbin{,}\kern-.1em}
\def\hk{\kern.12em\raise-1em\hbox{$\hat{\raise1em\hbox{,}}$}\kern.12em}
\newcommand{\initiate}{\setcounter{equation}{0}}
\begin{document}

\title{Metrics on the Real Quantum Plane}

\author{G. Fiore\footnote{Electronic
         mail:Gaetano.Fiore@na.infn.it}, $\strut^{1,2}$ 
         \, M. Maceda\footnote{Electronic
         mail:Marco.Maceda@th.u-psud.fr}, $\strut^{3}$ \,  
         J. Madore\footnote{Electronic
         mail:John.Madore@th.u-psud.fr}\hspace{3mm}$\strut^{3,4}$ \\[6pt]
        \and
        $\strut^1$Dip. di Matematica e Applicazioni, Fac.  di Ingegneria\\ 
        Universit\`a di Napoli, V. Claudio 21, 80125 Napoli
        \and
        $\strut^2$I.N.F.N., Sezione di Napoli\\
        Mostra d'Oltremare, Pad. 19, 80125 Napoli
        \and
        $\strut^3$Laboratoire de Physique Th\'eorique et Hautes Energies\\
        Universit\'e de Paris-Sud, B\^atiment 211, F-91405 Orsay
        \and
        $\strut^4$Humboldt Universit\"at zu Berlin, Institut f\"ur Physik\\
        Invalidenstrasse 110. D-10115 Berlin
        }

\date{}
\maketitle

\begin{abstract}
Using the frame formalism we 
determine some possible metrics and metric-compatible connections on the
noncommutative differential geometry of the real quantum plane. 
By definition a metric maps the tensor product of two
1-forms into a `function' on the quantum plane.
It is symmetric in a modified sense,
namely in the definition of symmetry one has to replace the permutator map
with a deformed map $\sigma$ fulfilling some suitable conditions. 
Correspondingly, also the 
definition of the hermitean conjugate of the tensor product of
two 1-forms is modified (but reduces to the standard one
if $\sigma$ coincides with the permutator).
The metric is real with respect to such modified
$*$-structure.
\end{abstract}

\vfill
\noindent
\noindent
\newpage

\initiate
\section{Introduction and notation}
\label{intro}

It is an old idea \cite{Sny47a, Sny47b} that a noncommutative
modification of the algebraic structure of space-time could provide a
regularization of the divergences of quantum field theory, because the
representations of noncommutative `spaces' have a lattice-like
structure. The main aim of noncommutative geometry \cite{Con94} is to endow
such an algebra with additional structures (starting from a
differential calculus), so as to build a bridge between
the algebra  and its `geometrical' interpretation.
Since there is no unique prescription how to do this, it
is useful to test possible prescriptions first on simpler models.

In this paper we choose as a noncommutative space(time) algebra
model the socalled real quantum or Manin plane \cite{Man88}, and as a 
differential calculus upon it the socalled Wess-Zumino calculus
\cite{WesZum90,PusWor89}. The latter is charaterized by the property
that the relations defining the module of 1-forms are covariant under the
action of the quantum group $SL_q(2)$ and homogeneous in
the generators. We adopt the noncommutative geometry formalism
of \cite{Mad89c,Mou95,DubMadMasMou95,DimMad96,DubMadMasMou96b}. 

We start with a brief description of the latter.
Let $\c{A}$ be an algebra with differential calculus
$\{\Omega^*(\c{A}),d\}$~\cite{Con94} (here $\Omega^*(\c{A})$
denotes the algebra of differential forms on $\c{A}$
and $d$ the exterior derivative acting on the latter)
and suppose that the calculus has a
{\it frame}~\cite{Mad89c,DimMad96}, i.e. a basis of 1-forms
$\theta^i$ ($i=1,2,...,n$)
which commute with the elements of the algebra,
\begin{equation}
\theta^i f = f \theta^i.                                       \label{3.1.51}
\end{equation}
The relation 
\begin{equation}
df =\theta^i e_i f                                            \label{3.1.52}
\end{equation}
(with $f\in\c{A}$) defines a set of derivations $e_i$ dual to
$\theta^i$,
from which it follows that the module structure of $\Omega^1(\c{A})$
is given by
$$
f dg = \theta^if e_i g , \qquad dg f = \theta^i (e_i g) f .
$$
We see that the $\c{A}$-bimodule $\Omega^1(\c{A})$ is free of rank~$n$
as a left or right module. It can therefore be identified with the
direct sum
\begin{equation}
\Omega^1(\c{A}) = \bigoplus_1^n \c{A}                           \label{3.1.53}
\end{equation}
of $n$ copies of $\c{A}$. In this representation $\theta^i$ is given
by the element of the direct sum with the unit in the $i$-th position
and zero elsewhere. We shall refer to the integer $n$ as the dimension
of the geometry. 

The wedge product $\pi$ in $\Omega^*(\c{A})$ fulfills
relations of the form
\begin{equation}
\theta^i \theta^j\equiv \pi(\theta^i \otimes_{\c{A}} \theta^j) =
P^{ij}{}_{kl} \theta^k \theta^l                   \label{3.1.62}
\end{equation}
(we omit the symbol $\wedge$ of the
wedge product), where $P$ is a projector
\begin{equation}
P^{ij}{}_{mn} P^{mn}{}_{kl} =  P^{ij}{}_{kl}              \label{3.1.61}
\end{equation}
with entries $P^{ij}{}_{kl} \in \c{Z}(\c{A})$.
If in particular the wedge product is such that
the $\theta^i$ anti-commute then $P$ is the antisymmetric projector
$$
P^{ij}{}_{kl} = 
{1\over 2} (\delta^i_k \delta^j_l - \delta^j_k \delta^i_l). 
$$
{F}rom (\ref{3.1.53}) it follows immediately that the algebra and its
differential calculus are related in a simple manner. Let $\exterior^*_P$
be the exterior algebra over $\b{C}^n$ with the wedge product
defined by~(\ref{3.1.62}). Then with the identification~(\ref{3.1.53})
it follows that one can write
\begin{equation}
\Omega^*(\c{A}) = \c{A} \otimes \exterior^*_P.                \label{rep}
\end{equation}
Since the exterior derivative of $\theta^i$ is a 2-form it can
necessarily be written as
$$
d\theta^i = 
- {1\over 2} C^i{}_{jk} \theta^j \theta^k.                     
$$
where, because of (\ref{3.1.62}), the structure elements
can be chosen to satisfy the constraints
$$
C^i{}_{jk} P^{jk}{}_{lm} = C^i{}_{lm}.
$$
It will also be convenient to introduce the quantities
\begin{equation}
C^{ij}{}_{kl} = \delta^i_k \delta^j_l - 2 P^{ij}{}_{kl}.       \label{3.1.73}
\end{equation}
Then from (\ref{3.1.61}) we find that 
\begin{equation}
C^{ij}{}_{kl}C^{kl}{}_{mn} = \delta^i_m \delta^j_n.            \label{3.1.74}
\end{equation}

For simplicity, we shall further assume that the $e_i$ are
inner derivations: $e_i f= [\lambda_i,f]$, $\lambda_i\in\c{A}$.  
{F}rom the $\theta^i$ we can construct a 1-form
\begin{equation}
\theta = - \lambda_i \theta^i                             \label{3.1.54}
\end{equation}
in $\Omega^1(\c{A})$ which plays the role of a Dirac
operator~\cite{Con94}:
\begin{equation}
df = - [\theta,f].                                       \label{3.1.55}
\end{equation}
One can show that
from the general consistency of the differential calculus it follows that
\begin{equation}
2P^{ij}{}_{kl} \lambda_i \lambda_j - F^i{}_{kl}\lambda_i - K_{kl} = 0   
                                                    \label{rellambda}
\end{equation}
for some array of elements $F^i{}_{jk},K_{kl}\in\c{Z}(\c{A})$.
In the cases which interest us here the latter vanish.

In order to consistently define a covariant
derivative we need to introduce~\cite{Mou95} a flip $\sigma$, i.e. a 
$\c{A}$-bilinear map
\begin{equation}
\Omega^1(\c{A}) \otimes_\c{A} \Omega^1(\c{A})
\buildrel \sigma \over \longrightarrow
\Omega^1(\c{A}) \otimes_\c{A} \Omega^1(\c{A}).                  \label{3.6.2} 
\end{equation}
In the case of the De-Rham calculus on the commutative 
algebra of functions on an oridinary manifold it reduces to
$\sigma(\omega\otimes_\c{A} \omega')=\omega'\otimes_\c{A} \omega$.
In terms of the frame it is given by
$S^{ij}{}_{kl} \in \c{Z}(\c{A})$ defined by
$$
\sigma(\theta^i \otimes_\c{A} \theta^j) = 
S^{ij}{}_{kl} \theta^k \otimes_\c{A} \theta^l.
$$
A covariant derivative on the module $\Omega^1(\c{A})$ is a map
\begin{equation}
\Omega^1(\c{A})
\buildrel D\over \longrightarrow
\Omega^1(\c{A}) \otimes_\c{A} \Omega^1(\c{A}).                  \label{defD} 
\end{equation}
satisfying
both a left and a right Leibniz rule. We use the ordinary left Leibniz
rule and define the right Leibniz rule as
\begin{equation}
D(\xi f) = \sigma (\xi \otimes_\c{A} df) + (D\xi) f           \label{3.6.3}
\end{equation}
for arbitrary $f \in \c{A}$ and $\xi \in \Omega^1(\c{A})$.
The connection 1-form 
$\omega^i{}_k\equiv\omega^i{}_{jk}\theta^j$ is defined by
\begin{equation}
D\theta^i=\omega^i{}_k\otimes_\c{A}\theta^k.
\end{equation}

We shall impose the condition 
\begin{equation}
\pi \circ (\sigma + 1) = 0                                     \label{3.6.5}
\end{equation}
that the antisymmetric part of a symmetric tensor vanish. This can be
considered as a condition on the product or on the flip.  In
ordinary geometry it is the definition of $\pi$; a 2-form can be
considered as an antisymmetric tensor. Because of this condition the
torsion is a bilinear map~\cite{DubMadMasMou96b}. The most general
solution can be written in the form
\begin{equation}
1 + \sigma = (1 - \pi) \circ \tau                  \label{s-t}
\end{equation}
where $\tau$ is an an arbitrary $\c{A}$-bilinear map.
Suppose that $\tau$ is invertible. Then
because of the identity
$$
1 = \pi + (1 + \sigma)\circ \tau^{-1}
$$
one can identify the second term on the right-hand side as the
projection onto the symmetric part of the tensor product.  The choice
$\tau = 2$ yields the value $\sigma = 1 - 2 \pi$.  If $\tau$ is not
invertible then there arises the possibility that part of the tensor
product is neither symmetric nor antisymmetric. Condition (\ref{3.6.5})
applied to the tensor product $\theta^i\otimes_\c{A} \theta^j$
becomes
\begin{equation}
P^{ij}{}_{lm}+S^{ij}{}_{hk}P^{hk}{}_{lm}=0.          \label{SP}
\end{equation}

If the flip is
such that in (\ref{rellambda}) $F^i{}_{jk}=K_{kl}=0$
one possible linear connection is
\begin{equation}
\omega^i{}_{jk} = 
\lambda_l (S^{il}{}_{jk} - \delta^l_j \delta^i_k).     \label{4.1.39}
\end{equation}
The corresponding connection 1-form is given by
\begin{equation}
\omega^i{}_k = \lambda_l S^{il}{}_{jk} \theta^j + \delta^i_k \theta.
                                                      \label{4.1.39'}
\end{equation}
The curvature of the covariant derivative $D$
defined in (\ref{4.1.39}) can be readily calculated. One finds the
expression
\begin{equation}
\frac 12 R^i_{jkl} = 
S^{im}{}_{rn} S^{np}{}_{sj} P^{rs}{}_{kl}\lambda_m \lambda_p.                           \label{curv}
\end{equation}
This can also be written in the form
$$
\frac 12 R^i_{jkl} = 
- S^{im}{}_{rn} S^{np}{}_{sj} S^{rs}{}_{uv} P^{uv}{}_{kl}
\lambda_m \lambda_p.
$$

In complete analogy with the commutative case a metric $g$ can be
defined as an $\c{A}$-bilinear, nondegenerate map~\cite{DubMadMasMou96b}
\begin{equation}
\Omega^1(\c{A}) \otimes_{\c{A}} \Omega^1(\c{A}) 
\buildrel g \over \longrightarrow \c{A}                        \label{3.4.1}
\end{equation}
and as such it can~\cite{CerHinMadWes99a} be used to define a
`distance' between `points'. It is important to notice here that the
bilinearity is an alternative way of expressing locality.  In ordinary
differential geometry if $\xi$ and $\eta$ are 1-forms then the value
of $g(\xi \otimes \eta)$ at a given point depends only on the values
of $\xi$ and $\eta$ at that point. Bilinearity, 
$$
g(f\xi \otimes_{\c{A}} \eta h)=f\,g(\xi \otimes_{\c{A}} \eta)\,h \qquad 
\qquad\forall f,h\in\c{A},
$$
is an exact expression
of this fact. In general the algebra introduces a certain amount of
non-locality via its nontrivial commutation relations and it is important to
assure that all geometric quantities be just that nonlocal and not
more.  Without the bilinearity condition it is not possible to
distinguish for example in ordinary space-time a metric which assigns
a function to a vector field in such a way that the value at a given
point depends only on the vector at that point from one which is some
sort of convolution over the entire manifold. 

We define frame components of the metric by
$$
g^{ij} = g(\theta^i \otimes_{\c{A}} \theta^j).
$$
They lie necessarily in the center $\c{Z}(\c{A})$ of the algebra.
The condition that (\ref{4.1.39}) be metric-compatible can be written as
\cite{DimMad96}
\begin{equation}
S^{im}{}_{ln} g^{np} S^{jk}{}_{mp} = g^{ij} \delta^k_l.      \label{3.6.40}
\end{equation}
As a a way to remember
this seemingly odd condition introduce a
`covariant derivative' $D_i X^j$ of a `vector' $X^j$.  The covariant
derivative $D_i (X^j Y)$ of the product of $X^j$ by a `field' $Y$ must
be then defined as
$$
D_i (X^j Y) = D_i X^j Y +S^{jl}{}_{im} X^m D_l Y
$$
since there is a `flip' as the index on the derivation crosses the
index on the first `vector'. If we apply again this rule to $Y=Y^k Z$,
with $Y^k$ also a `vector' and $Z$ another `field' we find
$$
D_i (X^j Y^k Z) = D_i (X^j Y^k) Z 
+ S^{jl}{}_{im} X^m Y^p S^{kn}{}_{lp} D_n Z.
$$
Since $g^{jk}$ is a `tensor', the `crossing rule` is the same as for
$X^jY^k$:
$$
D_i (g^{jk} Z) = (D_i g^{jk}) Z +S^{jl}{}_{im} g^{mp}S^{kn}{}_{lp} D_n Z.
$$
Therefore (\ref{3.6.40}) is equivalent to the usual condition,
$$
D_i g^{jk} = 0, 
$$
that the connection be compatible with the metric.

We shall require that the metric be symmetric in the sense
\begin{equation}
g \circ \pi = 0                                              \label{sym-cond}
\end{equation}
that it annihilates the 2-forms. 
This condition
applied to the tensor product $\theta^i\otimes_\c{A} \theta^j$
becomes
\begin{equation}
P^{ij}{}_{lm}g^{lm}=0.          \label{Pg}
\end{equation}

\bigskip
Let us now briefly summarize the additional conditions which
arise from the requirement of existence of $*$-structures.
Assume $\c{A}$ is a $*$-algebra.
If \cite{Wor89,Con95}
the $*$-structure of $\c{A}\equiv\Omega^0(\c{A})$ can be extended
to a $*$-structure of $\Omega^*(\c{A})$ and
\begin{equation}
(df)^*=df^*                                        \label{reald}
\end{equation}
the differential calculus is said to be real.
A sufficient condition for (\ref{reald}) to hold~\cite{FioMad98a}
is that the $\lambda_i$ are anti-hermitian (w.r.t. the $*$ of $\c{A}$)
and the $\theta^i$ are hermitean (w.r.t. the extension of $*$
to  $\Omega^*(\c{A})$), so that the `Dirac operator' $\theta$
is anti-hermitean. 

To obtain a real covariant derivative it is necessary first of all that
the flip $\sigma$ satisfies a reality constraint (see 
Ref.~\cite{FioMad98a}), which
takes the simple form
\begin{equation}
(S^{ji}{}_{kl})^* S^{lk}{}_{mn} = \delta^i_m \delta^j_n        \label{j-s}
\end{equation}
if $(\theta^i)^*=\theta^i$.  Moreover, the connection 1-form
$\omega^i_k$ and the flip $\sigma$ 
must satisfy a condition~\cite{FioMad98a} involving both, which we
do not report here because it
is automatically satisfied in the case of the connection (\ref{4.1.39}).
In order to define a real metric one has to
use $\sigma$ to impose the reality condition of Ref.~\cite{FioMad98a},
which takes the  simple form
\begin{equation}
S^{ij}{}_{kl} g^{kl} = (g^{ji})^*                             \label{her-f}
\end{equation}
in the case of a real frame. This is a
combination of a `twisted' symmetry condition and the ordinary
condition of reality on a complex matrix. It can also be written as an
ordinary condition of symmetry and a `twisted' definition of reality.
The map $\sigma$ is also involved~\cite{FioMad98a} in the reality
condition for the curvature or for the covariant derivative acting
on tensor powers of $\Omega^1(\c{A})$. The latter implies the former,
and takes the form of the braid equation,
\begin{equation}
S_{12}S_{23}S_{12}=S_{23}S_{12}S_{23}                     \label{braid}
\end{equation}
where
$$
(S_{12})^{abc}{}_{def} := S^{ab}{}_{de} \delta^c_f, \qquad 
(S_{23})^{abc}{}_{def} :=\delta^a_dS^{bc}{}_{ef}
$$

The `infinitesimal distance' $ds$ corresponding to the metric $g$
is introduced through the relation
\begin{equation}
ds^2=g_{ij}\theta^i\otimes_{\c{A}}\theta^j,           \label{defds}
\end{equation}
where $g_{ij}\in \c{A}$ are the matrix elements of the inverse matrix
of $\Vert g^{ij}\Vert$.
Every representation of $\c{A}$ yields a distance between `points'
because of~(\ref{rep}). Let $dt = \xi_i \theta^i \in \Omega^1(\c{A})$
be an exact form, which we can think of as an infinitesimal
displacement along an axis $t$ and suppose that $\ket{p}$ is a common
eigenvector of all the $\xi_i$: $\xi_i \ket{p} = \t{\xi}_i \ket{p}$.
This would be the case for example if only one of them is not equal to
zero.  We define the element of distance $\delta s$ along the
`coordinate' $t$ at the state $\ket{p}$ by the equation
$$
(\delta s)^2 = \bra{p}ds^2\ket{p} = g^{ij} \t{\xi}_i\t{\xi}_j.
$$
Let $\kbar$ be the length scale at which points become fuzzy and
$K^{-1}$ the scale at which the curvature effects become important.
The definition of $g$ which we have given is unambiguous but the
interpretation of the norm $|\delta s|^2$ of an infinitesimal displacement
as a distance can be only made within the range
$$
\kbar << |\delta s|^2 << K^{-1}.
$$
If the displacement is too small then the points are not defined;
if it is too large then an integral must be taken. The second problem
was solved by Leibniz/Newton; the first is a feature, not a bug, of
noncommutative geometry. We are especially interested in the region
$|\delta s|^2 \simeq \kbar$ where the noncommutative effects become of
interest.

There exist other definitions of distance. One
proposal~\cite{ConLot92,Lan97,FigGraVar00} uses the Dirac operator to
define distance on the space of pure states. Several
authors~\cite{FicLorWes96,LorWeiWes97b} do not consider the
bilinearity condition we have imposed as important and
several~\cite{PodWor90,Maj93,AscCas96,Pod96,AscCasSca99,Maj99,KosLukMas99}
consider the invariance under the coaction of a quantum group as
essential.

It is sometimes convenient to write the metric as a sum
$$
g^{ij} = g_S^{ij} +g_A^{ij}  
$$
of a symmetric and an antisymmetric part (in the usual sense of the
word) The inverse matrix we write as a sum
$$
g_{ij} = \eta_{ij} + B_{ij}
$$
of a symmetric and an antisymmetric term. We shall choose as
normalization when possible the condition that $\eta_{ij}$ be the
standard Minkowski or euclidean form.

\initiate
\section{The Wess-Zumino calculus}

The extended real quantum plane is the $*$-algebra $\c{A}$ generated by
hermitian elements $(x^i) = (x,y)$, 
\begin{equation}
x^*=x \qquad\qquad y^*=y,                              \label{real*}
\end{equation}
together with their inverses, fulfilling the relation
\begin{equation}
xy = \t{q}yx                                      \label{cr}
\end{equation}
with $|\t{q}|=1$ and $q\neq \pm 1$,
as well as the usual relations between inverses. We call it
{\it extended} because in the original version \cite{Man88}
the inverses $x^{-1},y^{-1}$ were not included; the
word {\it real} refers to the $*$-structure (\ref{real*}). The center
of $\c{A}$ is trivial, $\c{Z}(\c{A})=\b{C}$.
We now show how the Wess-Zumino calculus \cite{WesZum90} fits
in the scheme described in the previous section.
We define, for $\t{q}^4 \neq 1$,
$$
\lambda_1 = - \epsilon_1 {\t{q}^4 \over \t{q}^4 - 1} 
x^{-2} y^2,       \qquad
\lambda_2 = \epsilon_2 {\t{q}^2 \over \t{q}^4 - 1} 
x^{-2}.          
$$
There is an ambiguity in this definition due to the fact that the
defining relations (\ref{cr}) are homogeneous and which we reduce
to a sign: $\epsilon_a = \pm 1$. The extra minus is a `historical
convenience'.  The important fact is that the $\lambda_a$ are singular
in the limit $\t{q} \to 1$ and that they are anti-hermitian if $\t{q}$
is of unit modulus, as we are assuming.  We find for $\t{q}^2 \neq -1$
\begin{equation}
\begin{array}{ll}
e_1 x = \epsilon_1 \displaystyle{{\t{q}^2 \over (\t{q}^2 + 1)}}  
x^{-1} y^2,
&e_1 y =  \epsilon_1\displaystyle{{\t{q}^4 \over \t{q}^2 + 1} 
x^{-2}} y^3,\\[10pt]
e_2 x = 0, &e_2 y = 
- \displaystyle{\epsilon_2{\t{q}^2 \over \t{q}^2 + 1}} 
x^{-2} y.
\end{array}                                                    \label{4.1.50}
\end{equation}
These derivations are again extended to arbitrary polynomials in the
generators by the Leibniz rule. Using them and (\ref{3.1.54}),
(\ref{3.1.55}) we find
\begin{equation}
dx = {\t{q}^2 \over (\t{q}^2 + 1)} x^{-1} y^2 \epsilon_1\theta^1, \qquad
dy = {\t{q}^2 \over \t{q}^2 + 1} x^{-2} y 
(\t{q}^2 y^2 \epsilon_1 \theta^1 - \epsilon_2\theta^2)    \label{Dim}
\end{equation}
and solving for the $\theta^i$ we obtain
$$
\epsilon_1\theta^1 = (\t{q}^2 + 1) x y^{-2} dx,  \qquad
\epsilon_2\theta^2 = - (\t{q}^2 + 1) x (x y^{-1} dy - dx). 
$$
The module structure which follows from the condition
(\ref{3.1.51}) that the $\theta^i$ commute with the elements of the
algebra is equivalent to the Wess-Zumino relations~\cite{WesZum90}
\begin{equation}
\begin{array}{ll}
xdx = \t{q}^2 dx x, &x dy = 
\t{q} dy x + (\t{q}^2-1)dx y,   \\[4pt]
ydx = \t{q} dx y,   &y dy = 
\t{q}^2 dy y.
\end{array}                                                  \label{4.1.51}
\end{equation}
One can show that they are invariant under the coaction of the quantum
group $SL_q(2,\b{C})$.  This invariance was encoded in the choice of
$\lambda_a$.

Consider the elements 
\begin{equation}
u := \epsilon_2 \t{q}^{-2} x^{2}, \qquad 
v :=  \epsilon_1 x^2 y^{-2}. 
\end{equation}
We shall see that each of the four possible choices of sign
pairs corresponds to an identification of $x$ and $y$
as the coordinates of one of the four regions on $\b{R}^2$ defined by
the light cone of a metric with Minkowski signature. 
The $u, v$ fulfill the quadratic commutation relation 
\begin{equation}
uv=qvu
\end{equation}
where $q:=\t{q}^{-4}$. They and their inverses generate a slightly smaller
algebra than $\c{A}$.
One also finds that (\ref{4.1.51}) becomes
\begin{equation}
\begin{array}{ll}
u du = q^{-1} du u,      &u dv = q dv u,               \\[4pt]
v du = q^{-1} du v,      &v dv = q dv v. 
\end{array}                                            
\end{equation}
In terms of the new generators the $\theta^i$ become
\begin{equation}
\theta^1 =   q^{-1} vu^{-1} du,\qquad
\theta^2 =   u v^{-1} dv.               \label{thetauv}  
\end{equation}
What we have done in fact is use the $\lambda_a^{-1}$ as generators
of the algebra and the differential calculus; otherwise nothing has
been changed.  The form $\theta$ is most conveniently expressed in
terms of the $\lambda_a$. Since
\begin{equation}
\lambda_1 = \frac{1}{1-q^{-1}} v^{-1},                \qquad
\lambda_2 = -\frac 1{1-q^{-1}} u^{-1}                     \label{q-lambda}
\end{equation}
we find that
$$
\theta = {1 \over 1-q} (u^{-1}du-q v^{-1}dv).     
$$
It is an anti-hermitian closed form with vanishing square,
\begin{equation}
d \theta =0, \qquad \qquad( \theta)^2=0.          \label{dtheta=0}
\end{equation}
The volume element is a product of two exact forms:
$$
\theta^1 \theta^2 = du dv.
$$

The structure of the exterior algebra is given by the relations
\begin{equation}
(\theta^1)^2=0, \qquad 
(\theta^2)^2=0, \qquad
\theta^1\theta^2 + q\theta^2\theta^1 = 0.                \label{thth}      
\end{equation}
This can be written in the form (\ref{3.1.62}) with
\begin{equation}
P = \frac 12 \pmatrix{
 0  &    0      &  0  &  0  \cr
 0  &    1      & -q  &  0  \cr
 0  & -q^{-1}   &  1  &  0  \cr
 0  &    0      &  0  &  0
}.                                                          \label{P}
\end{equation}
If we reorder the
indices $(11,12,21,22) = (1,2,3,4)$ then the $C^{ij}{}_{kl}$
introduced in (\ref{3.1.73}) is given by the expression
$$
C = \pmatrix{
1 &0      &0 &0 \cr
0 &0      &q &0 \cr
0 &q^{-1} &0 &0 \cr
0 &0      &0 &1
}.                                                            
$$
That is, $C^{12}{}_{21} = q$ and $C^{21}{}_{12} = q^{-1}$.

The reality of the differential implies that the
structure elements must satisfy the conditions
$$
((C^i{}_{jk})^* + C^i{}_{jk}) P^{jk}{}_{lm} = 0
$$
from which it follows that
$$
(C^i{}_{21})^* = - C^i{}_{12} = q^{-1} C^i{}_{21}, \qquad
(C^i{}_{12})^* = - C^i{}_{21} = q C^i{}_{12}. 
$$
More precisely, the independent coefficients are given by
\begin{equation}
C^1{}_{12} = (q^{-1}-1) \lambda_2, \qquad 
C^2{}_{12} = (q^{-1}-1) \lambda_1.    \label{4.1.38}
\end{equation}
The $C^i{}_{jk}$ do not depend on the sign ambiguities.  With the
generators
\begin{equation}
t = \frac 1{\sqrt 2} (u + v), \qquad 
r = \frac 1{\sqrt 2} (u - v)              \label{tx-uv1}
\end{equation}
the four possible sign combinations can be written as 
$$
\epsilon_1 = \epsilon_2: \quad \sgn(t) =  \epsilon_1, \qquad
\epsilon_1 = -\epsilon_2:\quad \sgn(r) = \epsilon_2.
$$
We shall later in Section~5.1 introduce a light-cone and interpret
these relations in terms of space-like and time-like.

Introduce the notation
$$
X = \left(\begin{array}{c} t\\r \end{array}\right), \qquad
\Xi = \left(\begin{array}{c} dt\\dr \end{array}\right), \qquad
Q  = \left(\begin{array}{cc}  \cos (\pi\gamma) & i\sin (\pi\gamma) \\ 
    i\sin (\pi\gamma) & \cos (\pi\gamma) \end{array}\right)
\qquad q = e^{2\pi i\gamma}.
$$
Then $Q$ is unitary. The commutation relations in $\Omega^*(\c{A})$
can be written in the form
\begin{equation}
X^t (Q\sigma_2) X = 0, \qquad X \Xi^t = \Xi (Q^2 X)^t, \qquad 
\Xi^t Q \Xi = 0.                                                 \label{x-y}
\end{equation}
The $\sigma_2$ is the second Pauli matrix.

There are alternative $*$-structures which require a real $q$. One
can impose the conditions $u^*=v$, $v^*=u$. In terms of the original
variables $x$ and $y$ this implies that
$$
x^{*} = \pm {\t{q}}^{1/2} x y^{-1}, \qquad y^* = y.
$$
It follows that the frame satisfies
$$
(\theta^1)^* = \theta^2, \qquad (\theta^1)^* = \theta^2
$$
and so one can introduce a real frame by taking the real and
imaginary parts or consider the resulting structure as a $q$-deformed
complex line. This is better with the change of generators
\begin{equation}
t = \frac 1{\sqrt 2} (u + v), \qquad 
r = \frac i{\sqrt 2} (u - v).              \label{tx-uv1-bis}
\end{equation}
It is equivalent to a replacement $\gamma\mapsto i\gamma$ in the
formula (\ref{x-y}).

\initiate
\section{Representations}

An extensive discussion of the $*$-representations of the algebra $\c{A}$
for $|\t{q}|=1$ and $q\neq \pm 1$
has been given~\cite{Sch00}. We recall parts of it to illustrate our
interpretation of the geometry. It is easy to see that there can be no
normed basis with $u$ or $v$ diagonal.  Suppose in fact that there is
a basis with $v\ket{j} = v_j \ket{j}$. Since $v$ is hermitian the
eigenvalue $v_j\in\b{R}$. Using the commutation relations one sees
that $v(u\ket{j}) = q^{-1} v_j (u\ket{j})$ and so $u\ket{j}$ is also
an eigenvector with eigenvalue $q^{-1} v_j\notin\b{R}$.  One concludes
therefore that $u\ket{j}\notin\c{H}$. More specifically one can
consider $\c{H} = L^2(\b{R})$ with the plane-wave basis 
$\ket{k} = e^{ikx}$. The operator $u = - i\p_x$ is hermitian on a dense
subspace of $\c{H}$ and diagonal: $u\ket{k} = k \ket{k}$. We can
formally set
$$
v\ket{k} = \ket{q k} = e^{-iq kx}
$$
in order to have the correct commutation relations but $u$ is
not properly defined on the plane-wave basis.  

As solution to this problem we restrict our representation space to
the positive real line $\b{R}^+$ with free boundary condition at
$x=0$.  The Laplace transform replaces the Fourier transform and so we
choose as basis $\ket{k} = e^{-kx}$ for $k \in \b{C}$ with $\Re k > 0$.
We need in fact represent only one (at a time) of the four regions
defined by the light `cone' and we choose the one defined by 
$\epsilon_1 = \epsilon_2 = 1$. Our sign conventions were partly
dictated by the desire that this be the forward light-cone. We
choose~\cite{Sch00} then two positive real numbers $\alpha$ and
$\beta$ with $\alpha \beta = \gamma$ and we define on the Hilbert
space $L^2(\b{R}^+)$
$$
(uf)(x) = f(x+i\beta), \qquad (vf)(x) = e^{-2\pi\alpha x} f(x).
$$
Both $u$ and $v$ are formally hermitian and bounded. It is more
convenient to express them in terms of the Laplace transform,
which we recall is given by 
$$
F(k) = (Lf)(k) = \int_0^\infty f(x) e^{-k x}dx, \qquad
f(x) = (L^{-1}F)(x) = \frac 1{2\pi i}
\int_{a+i\infty}^{a-i\infty} F(k) e^{k x}dk
$$
where $a$ depends on the growth rate of the function. We have then
$$
(uF)(k) \equiv (L(uf))(k) = e^{i\beta k} F(k), \qquad
(vF)(k) \equiv (L(vf))(k) = F(k + 2\pi\alpha).
$$
In particular these transformation formulae are valid on the basis 
$\ket{k} = e^{-kx}$. The operators $u$ and $v$ are well-defined and
positive for $\Re k > 0$.

\section{The metrics and their connections}

We now determine some possible metrics and metric compatible 
connections on the
real quantum plane. We require them to fulfill all or at least
part of the conditions listed in section \ref{intro}, namely
(\ref{SP}), 
(\ref{braid}), 
(\ref{Pg}), 
(\ref{3.6.40}),
(\ref{j-s}), 
(\ref{her-f}). 

To shorten the notation we shall often
perform the following change of index notation:
$(11,12,21,22)\to(1,2,3,4)$.
Then the condition (\ref{3.6.40}) can be written
in the matrix form
\begin{equation}
\pmatrix{
S^1{}_1 & S^1{}_2 & S^1{}_3 & S^1{}_4  \cr
S^2{}_1 & S^2{}_2 & S^2{}_3 & S^2{}_4  \cr
S^3{}_1 & S^3{}_2 & S^3{}_3 & S^3{}_4  \cr
S^4{}_1 & S^4{}_2 & S^4{}_3 & S^4{}_4
}
\times \left(S_{(g)}\right) = 
\pmatrix{
g^1 &  0  & g^3 &  0  \cr
 0  & g^1 &  0  & g^3  \cr
g^2 &  0  & g^4 &  0  \cr
 0  & g^2 &  0  & g^4
}                                                          \label{4.1.40}
\end{equation}
where we have introduced the matrix $S_{(g)}$ defined by
\begin{equation}
S_{(g)} =
\pmatrix{
S^1{}_1 g^1 + S^1{}_2 g^3 & S^1{}_3 g^1 + S^1{}_4 g^3 & 
S^3{}_1 g^1 + S^3{}_2 g^3 & S^3{}_3 g^1 + S^3{}_4 g^3  \cr
S^1{}_1 g^2 + S^1{}_2 g^4 & S^1{}_3 g^2 + S^1{}_4 g^4 & 
S^3{}_1 g^2 + S^3{}_2 g^4 & S^3{}_3 g^2 + S^3{}_4 g^4  \cr
S^2{}_1 g^1 + S^2{}_2 g^3 & S^2{}_3 g^1 + S^2{}_4 g^3 & 
S^4{}_1 g^1 + S^4{}_2 g^3 & S^4{}_3 g^1 + S^4{}_4 g^3  \cr
S^2{}_1 g^2 + S^2{}_2 g^4 & S^2{}_3 g^2 + S^2{}_4 g^4 & 
S^4{}_1 g^2 + S^4{}_2 g^4 & S^4{}_3 g^2 + S^4{}_4 g^4 
}.
\end{equation}
Using the expression (\ref{P}) for $P$, the condition~(\ref{Pg}) 
becomes
\begin{equation}
g^2 = q g^3.                                                \label{eq:1}
\end{equation}
The consistency condition (\ref{3.6.5}) is equivalent to the
conditions
\begin{equation}
S^1{}_3 = q S^1{}_2, \qquad
S^2{}_3 = q(S^2{}_2 + 1), \qquad
S^3{}_3 = q S^3{}_2 - 1, \qquad
S^4{}_3 = q S^4{}_2.                                       \label{4.1.42}
\end{equation}

The equations to be solved then are Equations~(\ref{4.1.40}),
(\ref{eq:1}) and (\ref{4.1.42}). We are especially interested in real
solutions, which satisfy therefore also~(\ref{j-s}) and 
(\ref{her-f}). We have found
that there are several types of solutions, four of which we shall
describe in the following subsections. One can show that there are no
solutions with $\tau = 2$. A complete classification has
been given~\cite{Hie93} of the solutions to the braid equation as
well~\cite{GerGia97,AneArnChaDobMih00} as of those which satisfy a
weaker modified equation. 

If one considers locality as of importance only in the commutative
limit then there is no restriction on the coefficients of the metric,
except that they be local functions in this limit.  If one considers
locality as of importance even before the limit but is willing to
accept a metric which is real and symmetric only in the commutative
limit then the most general line element one can write is of the form
$$
ds^2 = g_{ij} \theta^i \otimes \theta^j.
$$
The $g_{ij}$ is a real symmetric matrix (in the sense we have defined it) 
and the moving frame $\theta^i$
is defined by
$$
\theta^1 = v u^{-1} du, \qquad \theta^2 =  u v^{-1} dv.
$$
The line element (\ref{defds}) becomes then
\begin{equation}
ds^2 = g_1 v^2 u^{-2} du^2 + 2 g_2 du dv + g_4  u^2 v^{-2} dv^2.\label{ds}
\end{equation}
The product here is the symmetrized tensor product; not the exterior
product. 

The associated metric connection is given by the structure functions
$$
C^1{}_{12} = u^{-1}, \qquad C^2{}_{12} = - v^{-1}.
$$

If we interpret the matrix $g_{ij}$ as the components of the Killing
metric on $SO(2)$ or $SO(1,1)$ 
then we can use it to calculate the
connection form. The result will be of the form
$$
\omega^i{}_j = A^i{}_{jk} u^{-1} \theta^k + B^i{}_{jk} v^{-1} \theta^k
$$
with $g_{ik} \omega^k{}_j$ antisymmetric in the two indices. 
The
Gaussian curvature $K$ is a second-order homogeneous polynomial in the
variables $u^{-1}$ and $v^{-1}$:
$$
K = \kappa_{11} u^{-2} + 2 \kappa_{12} u^{-1} v^{-1} + \kappa_{22} v^{-2}.
$$

\subsection{Solution I}
\label{SoluI}

A 1-parameter family of solutions of conditions (\ref{SP}), 
(\ref{Pg}), (\ref{3.6.40}),
can be found with a Minkowski-signature metric. For the
particular value $\zeta=0$ of the parameter also
the braid relation (\ref{braid}) and 
the reality conditions (\ref{j-s}), 
(\ref{her-f}) are fulfilled. These are the most interesting solutions.

With the convenient
normalization of the metric so that $g^3 = q^{-1/2}$ the flip is given
by the matrix
$$
S = \left(\begin{array}{cccc}
q & - q^{-1/2}\zeta  & - q^{1/2}\zeta 
               &   q^{-1}(q^2-1)^{-1}\zeta^2(q^2+1)  \\
0 &    0   & q & - q^{-1/2}\zeta          \\
0 & q^{-1} & 0 &   q^{-3/2}\zeta           \\
0 &    0   & 0 &   q^{-1} \end{array}\right),
$$
where $\zeta\in\b{C}$.
It tends to the ordinary flip as $q\to 1$ if $\zeta = 0$;
only for $\zeta = 0$ it is a
solution to the braid equation (\ref{braid}).  The corresponding
metric is given by
\begin{equation}
g^{ij} = \left(\begin{array}{cc}
(q-1)^{-1}\zeta & q^{1/2}   \\
    q^{-1/2}       &   0 
\end{array}\right).                                            \label{f-m-1}
\end{equation}
From~(\ref{eq:1}) one sees that it is $\sigma$-symmetric for all $g^1$
and real if $g^1 = 0$ (i.e. $\zeta = 0$). In this case $S$ is given by
\begin{equation}
S = \left(\begin{array}{cccc}
q &   0    & 0 & 0   \\
0 &   0    & q & 0   \\
0 & q^{-1} & 0 & 0   \\
0 &   0    & 0 & q^{-1}
\end{array}\right).                                          \label{4.1.43}
\end{equation}
The $\sigma$ and $\pi$ are related as in (\ref{s-t}) with 
$T^{ij}:=\tau(\theta^i\otimes_\c{A}\theta^j )$
\begin{equation}
T = \left(\begin{array}{cccc}
1+q & 0 & 0 & 0     \\
0   & 2 & 0 & 0     \\
0   & 0 & 2 & 0     \\
0   & 0 & 0 & 1+q^{-1}
\end{array}\right).                                         
\end{equation}
The fact that $T$ is not proportional to the identity is due to the
fact that the map $(1+\sigma)/2$ is not a projector and that we would
like it to act as such and be the complementary to $\pi$.  The metric
matrix is of indefinite signature and in `light-cone' coordinates.  If we use
the expression $q = e^{2\pi i\eta}$ we find that
\begin{equation}
g_S^{ij} = \cos (\pi\eta) \left(\begin{array}{cc}
   0 &  1   \\
   1 &  0 
\end{array}\right), \qquad 
g_A^{ij} = i \sin (\pi\eta) \left(\begin{array}{cc}
   0 &  1   \\
   -1 &  0 
\end{array}\right).
\end{equation}
The inverse metric components are defined by the equation
$$
g_{ij} g^{jk} = \delta^k_i.
$$
This matrix also can be split. If we rescale so that the symmetric
part is of the standard form we find
$$
(\eta_{ij} )= \left(\begin{array}{cc}
   0 &  1   \\
   1 &  0 
\end{array}\right), \qquad 
(B_{ij} ) = i \tan (\pi\zeta) \left(\begin{array}{cc}
   0 &  1   \\
  -1 &  0 
\end{array}\right).
$$
For the choice (\ref{4.1.43}) of the flip (i.e. for $\zeta=0$)
the metric connection~(\ref{4.1.39'}) is given by
$$
(\omega^i{}_j) = (1-q)\left(\begin{array}{cc}
1 &0 \\[4pt] 0 & - q^{-1}
\end{array}\right) \theta,
$$
and has vanishing curvature, 
because of the identities (\ref{dtheta=0}) and (\ref{braid}). This can be
shown by an argument already used in \cite{FioMad99,CerFioMad00a}.
In other words, in this case the quantum plane is flat.
In the commutative limit the line element is given by
$$
ds^2 = g_{ij}\theta^i \otimes \theta^j = 
2 \theta^1 \otimes\theta^2 = 
2 du \otimes dv = dt^2 - dr^2.
$$
The frame is singular along the light cone through the origin
[see (\ref{thetauv})].
Suppose $\epsilon_1 = \epsilon_2 =1$. If in a representation one
forces $x$ and $y$ to be hermitian then the $u$ and $v$ must be
positive operators. One concludes then that $t > |r|$; the geometry
describes only the forward light-cone through the origin. The other
three regions are given by the other three possible combinations of
signs.

\subsection{Solution II}

A family of solutions defined by flips which
are solutions to
(\ref{SP}), (\ref{j-s}), 
but not to the braid equation (\ref{braid}) is given by
\begin{equation}
S = \left(\begin{array}{cccc}
-q^2 &    0    &      0     & 0\\[2pt]
  0  &    0    &      q     & 0\\[2pt]
  0  & -q^{-2} &  -1-q^{-1} & 0\\[2pt]
  0  &    0    &      0     & q^{-1}
\end{array}\right)
\end{equation}
The metric is given again by~(\ref{f-m-1}) with $\zeta=0$, and
fulfills (\ref{3.6.40}), (\ref{Pg}), but not (\ref{her-f}). 
The metric connection~(\ref{4.1.39'}) is
$$
(\omega^i{}_j) = (1+q^2)\left(\begin{array}{cc}
1 & 0\\ 0 & q^{-2} \end{array}\right) \theta
+ (1+q^{-1}) \left(\begin{array}{cc} 
0&   0 \\[4pt] 
-1&  0 \end{array}\right)\lambda_1\theta^2
+ (q+1)\left(\begin{array}{cc} 
q & 0 \\[4pt] 
0 & q^{-2} \end{array}\right)\lambda_2\theta^2.
$$
The curvature Curv is equal to
$$
\Omega^i{}_j = -(q^2-1)q^{-3}(1+q+q^2)
\left(\begin{array}{cc} 0&0\\1&0\end{array}\right)
(\lambda_1)^2 \theta^1\theta^2.
$$
It diverges as $(q-1)^{-1}$ when $q\to 1$. This is then the case of a
regular metric which has a singular metric connection.

\subsection{Solution III}

A third family,
\begin{equation}
S = {1 \over q^2 + 1}\pmatrix{
2 q     & 0       & 0       & 1 - q^2   \cr
0       & 1 - q^2 & 2 q     & 0         \cr
0       & 2 q     & q^2 - 1 & 0         \cr
q^2 - 1 & 0       & 0       & 2 q
},                                                            \label{4.1.44}
\end{equation}
$$
g^{ij} = \left(\begin{array}{cc}1 & 0 \\[2pt] 0 & 1
\end{array}\right),
$$
fulfills
(\ref{SP}), (\ref{Pg}), (\ref{3.6.40}), the reality condition
(\ref{j-s})  but not the one (\ref{her-f}) nor the braid relation
(\ref{braid}). The latter are fulfilled
 for $q = \pm 1$. For $q = -1$ this means the connection form is
imaginary in the usual sense of the word (since so are the $\lambda_i$). 

The compatible connection (\ref{4.1.39'}) form is 
$$
\displaystyle{
(\omega^i{}_j) = \frac{(q-1)^2}{q^2+1}
 \delta^i_j \theta + \frac{q^2-1}{q^2+1}
\left(\begin{array}{cc}0 & -1 \\[2pt] 1 & 0 \end{array}\right) 
(\lambda_2 \theta^1 + \lambda_1 \theta^2).
}
$$
The curvature 2-form is
$$
\displaystyle{
(\Omega^i{}_j) = \frac{(q^2-1)}{(q^2+1)^2}
\left\{- q^{-1}(q^2-1)^2 \delta^i_j
\lambda_1\lambda_2  + 2 (q - 1) 
\left(\begin{array}{cc}0 & -1 \\[2pt] 1 & 0
\end{array}\right)\Big((\lambda_1)^2 + (\lambda_2)^2\Big)\right\} 
\theta^1 \theta^2.
}
$$
In the limit $q\to 1$ this becomes
$$
\displaystyle{
(\Omega^i{}_j )=  
\left(\begin{array}{cc}0 & -1 \\[2pt] 1 & 0
\end{array}\right)(u^{-2} + v^{-2}) \theta^1 \theta^2.
}
$$

\subsection{The $\h{R}$-matrix `solution'}

Finally one might ask whether one can find a solution $(S,g)$ using
the formalism of Faddeev {et al.}~\cite{FadResTak89}, as has been
done~\cite{FioMad99,CerFioMad00a} for the $q$-euclidean `spaces´
$\b{R}^n_q$ with $n>2$.  This would imply an $S$ proportional to the
braid matrix $\h{R}$ of $SL_q(2)$ or to its inverse. One can show that
there is a solution only if one admits non-symmetry metrics.

We recall that the braid matrix which defines the Hopf algebra
$SL_q(2)$
$$
\h{R}_q=\left(\begin{array}{cccc}
 q &     0    & 0 & 0   \\
 0 & q-q^{-1} & 1 & 0   \\
 0 &     1    & 0 & 0   \\
 0 &     0    & 0 & q
\end{array}\right)
$$
fulfills the braid relation, admits the projector decomposition
$$
\h{R}_q=qP_{s,q}-q^{-1}P_{a,q}
$$
and fulfills the (\ref{braid}) relations
\begin{equation}
\hat R_q^{\pm 1}{}^{ij}{}_{hk}\varepsilon_q^{kl}
\hat R_q^{\pm 1}{}^{rs}{}_{jl}=
q^{\mp 1}\varepsilon_q^{ir}\delta^s_h,
\qquad\qquad\hat R_q^{\pm 1}{}^{ij}{}_{hk}\varepsilon_q^{hk}=
-q^{\mp 1}\varepsilon_q^{ij},                          \label{rrr}
\end{equation}
where $\varepsilon_q^{ij}$ is the $q$-deformed epsilon tensor 
$$
\varepsilon_q^{ij} = \left(
\begin{array}{ll}
     0   & -q^{-1/2}\\
q^{1/2} &     0
\end{array}\right).
$$
So one finds
$$
P_{a,q}{}^{ij}{}_{hk} = (\varepsilon^{lm}\varepsilon_{lm})^{-1}
(\varepsilon^{ij}\varepsilon_{hk})=   
\frac{1}{q+q^{-1}}\left(\begin{array}{cccc}
0 &   0    &   0   &   0  \\[4pt]
0 & q^{-1} &  -1   &   0  \\[4pt]
0 &  -1    &   q   &   0  \\[4pt]
0 &   0    &   0   &   0
\end{array}\right).
$$

By a straightforward computation one can check that (\ref{thth}) 
can be given the form (\ref{3.1.62}) by setting
$$
P= P_{a,q^{-1}}.
$$
The first relation in (\ref{rrr}) suggests that we make the Ansatz
$S\propto \hat R_{q^{-1}}^{\pm 1}$, 
$g^{ij}\propto \varepsilon_{q^{-1}}^{ij}$, so that we can fulfill
(\ref{3.6.40}) at least up to a conformal factor.
Equation~(\ref{3.6.5}) fixes the first proportionality constant to be
either
$$
S=q^{-1}\hat R_{q^{-1}} \quad \mbox{or} \quad S=q(\hat R_{q^{-1}})^{-1}
$$
which respectively imply that
\begin{equation}
S^{im}{}_{ln} g^{np} S^{jk}{}_{mp} = q^{-1}g^{ij} \delta^k_l   \qquad
S^{im}{}_{ln} g^{np} S^{jk}{}_{mp} = qg^{ij} \delta^k_l ,  
\end{equation}
i.e. we indeed  fulfill
(\ref{3.6.40}) only up to a conformal factor $q^{\pm 1}$, and
\begin{equation}
S^{ij}{}_{hk}g^{hk}=-g^{ij}.                         \label{antis}
\end{equation}
This `antisymmetry' relation is to be contrasted with
Equation~(\ref{sym-cond}), which, with the above choice of $S$,
amounts to replacing at the rhs of (\ref{antis}) $-1$ respectively by
$q^{-2}$ or $q^2$, as can be seen writing $P$ as a combination of $S$
and of the identity matrix.  Using the fact that $|q|=1$ and $\hat
R_{q^{-1}}{}^{ij}{}_{hk}=\hat R_q^{-1}{}^{ji}{}_{kh}$ 
\cite{FadResTak89} one can easily
see that the reality conditions (\ref{j-s}) and (\ref{her-f}) are
satisfied.  The curvature (\ref{curv}) can be easily calculated to be
zero using the conditions $K^{ij}=0$ and $F^h_{ij}=0$ as well as the
fact that $P_q$ is a polynomial in $S$, which it turn fulfils the
braid equation.


\subsection{Other `solutions'}

There are a certain number of partial solutions which are
unsatisfactory for some reason or other. As an example, to underline
the possibility of exotic metrics which are both symmetric and
anti-symmetric according to our definitions, we consider $\sigma$
defined by the matrix
$$
S =   \left(\begin{array}{cccc}
     0      &   0  &  0   & \zeta \\[4pt]
     0      &  -1  &  0   &   0    \\[4pt]
     0      &   0  &  -1  &   0    \\[4pt]
\zeta^{-1}  &   0  &   0  &   0
\end{array}\right)
$$
where $\zeta \in \b{R}$ is a parameter.  This value of $S$ is a
solution to the braid equation. The $\sigma$ and $\pi$ are related as
in (\ref{s-t}) with (using the same conventions)
\begin{equation}
{\bf 1}+S=T = \left(\begin{array}{cccc}
      1     & 0 & 0 & \zeta     \\
      0     & 0 & 0 & 0     \\
      0     & 0 & 0 & 0     \\
\zeta^{-1}  & 0 & 0 & 1
\end{array}\right).                                         
\end{equation}
This means that $\tau$ is not invertible and the case is degenerate.
The unpleasent thing here is that $(1+\sigma)/2$ and $\pi$ do not add up to
the identity map.  The metric is given by
\begin{equation}
g^{ij} = i \left(\begin{array}{cc}
1   &  0   \\
0 & -\zeta^{-1} 
\end{array}\right).
\end{equation}
One has $\tau = 1 + \sigma$ and the flip is degenerate. Instead of
interchanging $g^2$ and $g^3$ as does the ordinary flip, it
interchanges $g^1$ and $g^4$. It also changes the sign, which accounts
for the $i$ in the metric components. Also $g \circ (1 + \sigma) = 0$
so in a certain sense the metric has vanishing symmetric as well as
antisymmetric parts. We refer to $\sigma$ nonetheless as a `flip'
because it satisfies~(\ref{3.6.5}).

The linear connection~(\ref{4.1.39}) is given by
$$
\omega^i{}_j = \delta^i_j \theta + 
\left(\begin{array}{cc} 0 &1 \\[4pt] 
-\zeta^{-1} & 0 \end{array}\right) 
(\zeta\lambda_1 \theta^2 - \lambda_2 \theta^1)
$$
The curvature is given by
$$
\Omega^i{}_j = q^{-1}(q^2-1) \delta^i_j \lambda_1\lambda_2 \theta^1 \theta^2
$$
The connection is singular in the commutative limit as is the
curvature.  Because of~(\ref{sym-cond}) it cannot be satisfied for any
curvature which is proportional to the metric.

\initiate
\section{Jordanian deformation}

It has been shown recently (See, for example, Aneva
\etal~\cite{AneArnChaDobMih00}) that the jordanian deformation is a
singular limit of a family of $q$ deformations.  The transformation
from the set of generators of one algebra to the other has also been
studied in some detail~\cite{Cha00}. We can now discuss to what extent
the limit can be understood in a geometric manner. We recall that the
jordanian deformation is defined using a parameter $h$ and that the
generators $(x^\prime, y^\prime)$ satisfy the commutation relations
$[x^\prime,y^\prime] = h y^{\prime 2}$. The differential calculus is
given by two elements $\lambda^\prime_a$ similar to the $\lambda_a$
which satisfy the $SL(2,\b{R})$ relation 
$[\lambda^\prime_1, \lambda^\prime_2] = \lambda^\prime_1$, a relation
which is not quadratic.  This must be compared with the quadratic
relation $\lambda_1 \lambda_2 = q^{-1}\lambda_2 \lambda_1$ satisfied
by the elements~(\ref{q-lambda}). We must find a smooth map from one
algebra into the other, that is, one which respects the commutation
relations between the elements which define the derivations dual to
the frame.  Consider~\cite{Cha00} the map
\begin{equation}
\lambda_1^\prime = h_0^{-1} \lambda_1, \qquad  
\lambda_2^\prime = h_0^{-1} \lambda_2 - \frac 12 h^{-1} h_0
\qquad h_0 = \frac{2h}{1-q}.                                  \label{Chak}
\end{equation}
This change defines a deformation of the differential calculus.  From
the commutation relations of the $\lambda_i$ we deduce that
$$
[\lambda_1^\prime, \lambda_2^\prime] = h_0^{-2}[\lambda_1,\lambda_2] = 
h_0^{-2}(1-q) \lambda_1 \lambda_2 = \lambda^\prime_1 + 
(1-q) \lambda_1^\prime \lambda_2^\prime.
$$
In the (singular) limit when $q \to 1$ the differential calculus
tends to that of the jordanian deformation.

The relations between the two calculi can be written in terms of a diagram
\begin{equation}
\begin{array}{ccc}
(x,\; y) &\longrightarrow
&(u,\;v) = (\epsilon_2q^{1/2} x^{2}, \; 
\epsilon_1  x^{2} y^{-2}).\vspace{8pt}\\
\downarrow &&\downarrow \vspace{8pt}\\ 
(x^\prime,\, y^\prime) &\longrightarrow
&(u^\prime,\; v^\prime) = 
(x^\prime y^{\prime -1} + \frac 12 h,\; y^{\prime -2}) \\
\end{array}                                                      \label{1.9}
\end{equation}
The two horizontal arrows are changes of generators. The two vertical
ones define a map between the two deformations.  In terms of the
generators $u$ and $v$ and their analogues~\cite{ChoMadPar98}
$u^\prime$ and $v^\prime$ for the jordanian deformation, the
map~(\ref{Chak}) can be written as
$$
u^\prime =  q u^{-1} - h_0, \qquad
v^\prime = -  q v^{-1}
$$
with $h_0 \to\infty$. It has been shown~\cite{ChoMadPar98} that the
local metric on the jordanian deformation is that of Lobachevsky. This
must be a limit of one of the family of metrics~(\ref{ds}).  The
Lobachevsky metric can be described with the line element 
$ds^{\prime 2} = v^{\prime -2} (du^{\prime 2} + dv^{\prime 2})$.  To
compare we write~(\ref{ds}) in the primed variables:
$$
ds^2 =  (u^\prime + h_0)^{-2} v^{\prime -2} [q^2\,g_1du^{\prime 2}
- 2 g_2 du^\prime dv^{\prime}
+ q^{-2}g_4  dv^{\prime 2}].
$$
We see than that we must choose $g_2 = 0$ and let $g_1, g_4 \to\infty$
with the constraint
$$
g_1 h^{-2}_0 = g_4 h^{-2}_0 = 1.
$$
The quantum-plane metric belongs to the family III. Another
interesting metric obtained in the same limit is with $g_1 = g_4 = 0$
and $g_2\to\infty$ so that $g_2 h^{-2}_0 = 1$:
$$
ds^2 = - 2 v^{\prime -2} du^\prime dv^{\prime} = - 2 du^\prime dv.
$$
This solution belongs to the family I.

\initiate
\section{Patching}

Let us consider now the solutions {\it I} found in section
\ref{SoluI}.
To each of the four regions defined by the light cone through the
origin in two dimensions we have associated an algebra, a 
differential calculus and a metric, but none is
complete as `manifold'.  
From the form of the metric we see that this can be done
using the generators $(t,r)$ or $(u,v)$ but that the generators
$(x,y)$ are singular on the cone. 

The patching is done~\cite{Sch00} by
extending the domain of definition of $u$ for example to negative
eigenvalues. The frame $\theta^i$ is also singular on the cone but the
equivalent frame $du^i$ is quite regular. We can write
$\theta^i = \Lambda^i{}_j du^j$ where
$$
\Lambda^i{}_j = \sqrt q \left(
\begin{array}{cc} vu^{-1} & 0 \\ 0 & uv^{-1}
\end{array}\right)
$$
is a local Lorentz transformation in the commutative limit. 

\initiate
\section{Discussion}

We have given a partial classification of the solutions to the three
conditions of metric compatibility~(\ref{3.6.40}),
symmetry~(\ref{sym-cond}) and the consistency condition~(\ref{3.6.5}),
as well as the reality conditions (\ref{j-s}), 
(\ref{her-f}), and the braid relation
(\ref{braid}),
without due regard to quantum covariance.  In fact we could show that
there was no solution which respected a coaction of the quantum group.
A similar problem was found by Cotta-Ramusino \& Rinaldi in trying to
construct holonomy groups~\cite{CotRin92}.  Written in terms of the
components in the frame basis one sees that $S^{ij}{}_{kl}$ has 16
unknowns and $g^{ij}$ has 4 unknowns. The condition~(\ref{3.6.5})
gives 4 equations and metric compatibility gives 16 equations.  So a
naive computation would say that the solution is unique up to a
rescaling of $g^{ij}$, which is not fixed by the equation. We have
indeed found a finite set of solutions.

Another conclusion concerns the uniqueness of the vacuum. It has been
claimed~\cite{Mad99c} that within the context of the present formalism
there is essentially a unique differential calculus which has
associated to it a given metric, unique that is up to a choice of norm
on the frame.  This statement needs qualification since we have here
shown that the quantum plane is naturally endowed with the
Lorentz-signature flat metric and it is known that the same is true of
the Heisenberg algebra with its natural differential calculus.

\initiate
\section*{Acknowledgment} The authors would like to thank
A. Chakrabarti for enlightening conversations.  One of them (JM) 
was supported by the Deutsche Forschungsgemeinschaft and he would
like to thank Dieter L\"ust for his hospitality at the Institut f\"ur
Physik, Berlin, were part of this research was carried out. Also MM
would like to thank the National Council of Science and Technology
(CONACyT, M\'exico) for financial support.

\setlength{\parskip}{5pt}

\providecommand{\href}[2]{#2}\begingroup\raggedright\endgroup


\end{document}